\documentclass[reqno]{amsart}

\usepackage{amsopn} 
\usepackage{amssymb} 
\usepackage{xypic} 
\usepackage{graphicx} 
\usepackage[mathcal]{euscript} 

\numberwithin{equation}{subsection}
\newtheorem{introThm}{Theorem}
\newtheorem{thm}[equation]{Theorem}
\newtheorem{lem}[equation]{Lemma}
\newtheorem{prop}[equation]{Proposition}
\theoremstyle{remark}
\newtheorem*{remark}{Remark}
\newtheorem{convention}[equation]{Convention}
\newtheorem*{remarks}{Remarks}

\newcommand{\bk}{\Bbbk}
\newcommand{\sce}{\mathcal{E}}
\newcommand{\scg}{\mathcal{G}}
\newcommand{\g}{\textnormal{g}}
\newcommand{\mfa}{\mathfrak{a}}
\newcommand{\mfb}{\mathfrak{b}}
\newcommand{\mfm}{\mathfrak{m}}
\newcommand{\mfn}{\mathfrak{n}}
\newcommand{\bbz}{\mathbb{Z}}
\newcommand{\bbn}{\mathbb{N}}
\newcommand{\dev}{\varepsilon}
\DeclareMathOperator{\upP}{\tn{P}}

\DeclareMathOperator{\upH}{H}
\DeclareMathOperator{\upT}{T}
\DeclareMathOperator{\upU}{U}
\DeclareMathOperator{\ld}{ld}
\DeclareMathOperator{\lin}{lin}
\DeclareMathOperator{\reg}{reg}
\DeclareMathOperator{\codepth}{codepth}
\DeclareMathOperator{\depth}{depth}
\DeclareMathOperator{\polreg}{pol \ reg}
\newcommand{\surj}{\twoheadrightarrow}
\newcommand{\inj}{\hookrightarrow}
\newcommand{\ps}[1]{{}_{#1}}
\newcommand{\tn}[1]{\textnormal{#1}}
\newcommand{\cl}[1]{\overline{#1}}
\newcommand{\ang}[1]{\left\langle #1 \right\rangle}
\newcommand*{\bt}{{\scalebox{0.6}{$\bullet$}}}
\newcommand{\bd}{\partial}
\DeclareMathOperator{\gr}{gr}
\DeclareMathOperator{\codim}{codim}

\DeclareMathOperator{\im}{im}

\DeclareMathOperator{\Ext}{Ext}
\DeclareMathOperator{\grExt}{grExt}
\DeclareMathOperator{\grHom}{grHom}
\DeclareMathOperator{\grSoc}{grSoc}
\DeclareMathOperator{\Hom}{Hom}
\DeclareMathOperator{\edim}{edim}

\title[Minimal multiplicity]{Homological criteria for minimal multiplicity}
\author{John Myers}
\address{Department of Mathematics, SUNY Oswego, Oswego, New York, USA}
\email{john.myers@oswego.edu}

\begin{document}

\begin{abstract}

Lower bounds on Hilbert-Samuel multiplicity are known for several types of commutative noetherian local rings, and rings with multiplicities which achieve these lower bounds are said to have minimal multiplicity. The first part of this paper gives characterizations of rings of minimal multiplicity in terms of the Ext-algebra of $R$; in particular, we show that minimal multiplicity can be detected via an Ext-algebra which is Gorenstein or Koszul AS-regular. The second part of this paper characterizes rings of minimal multiplicity via a numerical homological invariant introduced by J. Herzog and S. B. Iyengar called linearity defect. Our characterizations allow us to answer in two special cases a question raised by Herzog and Iyengar.

\end{abstract}
\maketitle
\thispagestyle{empty}

This paper is concerned with three lower bounds on the Hilbert-Samuel multiplicities of certain commutative noetherian local rings: one bound for the class of Cohen-Macaulay rings, one for Gorenstein rings, and one for complete intersections. Rings whose multiplicities achieve these lower bounds (i.e., rings of \textit{minimal multiplicity}) are known to have several desirable properties --- the ones that concern us are homological in nature, and the main contribution of this paper is to show that these properties provide sufficient conditions for a ring to have minimal multiplicity.

After some preliminary results are described in Section 1, our main results begin in Section 2 with the $\Ext$-algebra $\sce_R$ of a commutative noetherian local ring $R$, a noncommutative graded algebra constructed from $R$ via homological algebra. We prove

\begin{introThm}\label{thm:1}
Let $R$ be a commutative artinian local ring of embedding dimension $\geq 2$. The following are equivalent:
\begin{enumerate}
\item The ring $R$ is Gorenstein of minimal multiplicity.
\item The ring $R$ is Gorenstein, the algebra $\sce_R$ has global dimension $2$, and there is an isomorphism $\sce_R \cong (R^\square)^!$ of connected algebras.
\item The algebra $\sce_R$ is quadratic Gorenstein of global dimension $2$.
\end{enumerate}
\end{introThm}

The notation $R^\square$ denotes the quadratic part of the tangent cone of $R$, and $(-)^!$ denotes the quadratic (or Koszul) dual of a quadratic algebra. The isomorphism in (2) is the essential link which allows us to pass information from an $\Ext$-algebra to the local ring from which it was built; this isomorphism follows from computations made by Sj\"odin in \cite{Sjodin1976}. We frame the discussion leading up to Theorem \ref{thm:1} in terms of this isomorphism, asking under what conditions on $\sce_R$ and $R$ does it exist.

The Gorenstein condition in statement (3) is a noncommutative analog of the same condition in commutative algebra. The proof of Theorem \ref{thm:1} flows through general arguments and results which have nothing to do with commutative local rings; for example, the proof uses the fact Gorenstein algebras are quadratic dual to (graded) Frobenius algebras.

With artinian Gorenstein rings addressed, Section 2 moves on to artinian complete intersections. The analog to Theorem \ref{thm:1} is

\begin{introThm}\label{thm:2}
Let $R$ be a commutative artinian local ring and set $m=\edim{R}$. The following are equivalent:
\begin{enumerate}
\item The ring $R$ is a complete intersection of minimal multiplicity.
\item The ring $R$ is a complete intersection, the algebra $\sce_R$ has global dimension $m$, and there is an isomorphism $\sce_R \cong (R^\square)^!$ of connected algebras.
\item The algebra $\sce_R$ is Koszul AS-regular of global dimension $m$.
\end{enumerate}
\end{introThm}

The AS-regular algebras in (3) are from noncommutative algebraic geometry, introduced and studied by Artin and Schelter (hence the name) in \cite{ArtinSchelter1987} as noncommutative analogs of commutative polynomial algebras. These algebras are Gorenstein, but also have polynomial growth and finite global dimension. It has been known since Gulliksen showed in \cite{Gulliksen1980} that the $\Ext$-algebras of all complete intersections (non necessarily artinian) have polynomial growth, so our main contribution and the bulk of the proof of Theorem 2 concerns the Gorenstein condition. Our proof makes use of an important extra feature of $\sce_R$, namely that it is the universal enveloping algebra of a graded Lie algebra.

Section 3 begins with the definition of a numerical homological invariant of a finitely-generated module $M$ over a commutative noetherian local ring $R$ introduced by Herzog and Iyengar in \cite{HerzogIyengar2005} and called the linearity defect of $M$. This invariant is either a nonnegative integer or $\infty$, and those modules with finite linearity defect are described by Herzog and Iyengar to be ``modules in whose minimal free resolution the linear part predominates.'' The linearity defect can also be defined for commutative graded connected algebras over a field $k$, and here it was noted by Herzog and Iyengar that the linearity defect of $k$ has a rigidity property: If this linearity defect is finite, then it must be zero. Whether or not this is also true in the local case is an open question, with substantial progress being made by \c Sega in \cite{Sega2013} who uncovered a link with minimal multiplicity. \c Sega's results depend on the tangent cone of the local ring under consideration being Cohen-Macaulay, but by using a result of Levin \cite{Levin1985} we show that this Cohen-Macaulay hypothesis is not necessary and prove

\begin{introThm}\label{thm:3}
Let $R$ be a commutative noetherian local ring with residue field $k$. If either $R$ is a complete intersection or is both Cohen-Macaulay and Golod, then the following are equivalent:
\begin{enumerate}
\item The ring $R$ is Koszul (i.e., the $R$-module $k$ has linearity defect zero).
\item The ring $R$ is Fr\"oberg.
\item The ring $R$ has minimal multiplicity.
\item The $R$-module $k$ has finite linearity defect.
\end{enumerate}
\end{introThm}

The Fr\"oberg condition in (2) expresses a relation between the Hilbert series of $R$ and the Poincar\'e series of the residue field $k$. It is known that all Koszul local rings are Fr\"oberg, but it is an open question as to whether or not the converse is true in general. It is known for complete intersections, but to the author's knowledge the equivalence (1) $\Leftrightarrow$ (2) in Theorem \ref{thm:3} for rings which are Cohen-Macaulay and Golod is new.

\section{Definitions and preliminaries}

This section fixes notation and terminology, and states a few results which will be used later. Our three main references are \cite[Chapter I]{BellamyEtAl2016} and \cite[Chapters 1 and 2]{PolishchukPositselski2005} for graded algebra, and \cite{BrunsHerzog1998} for commutative algebra.

\subsection{Connected algebras and graded modules}

A \textit{connected} algebra over a field $\bk$ is an $\bbn$-graded $k$-algebra $\bigoplus_{j\geq 0} A_j$ such that each $A_j$ is finite-dimensional over $\bk$ (i.e., $A$ is \textit{locally finite-dimensional}) and the structure map $\bk \to A$ induces an isomorphism $A_0 \cong \bk$ of vector spaces. A \textit{graded module} over a connected algebra is a $\bbz$-graded module $M=\bigoplus_{j\in \bbz} M_j$ which is locally finite-dimensional and $M_j=0$ for all $j\ll 0$ (i.e., $M$ is \textit{bounded below}). We say $M$ has \textit{polynomial growth} if there is a positive constant $a$ and a positive integer $b$ such that $\dim_k M_j \leq a j^b$ for all $j>0$. We put
	\[\sup{M} = \sup \{ j\in \bbz : M_j \neq 0\},
	\]
with the convention that $\sup{M}=+\infty$ if the set $\{j\in \bbz: M_j \neq 0\}$ is unbounded above. If $n$ is an integer, the \textit{$n$-th internal shift} of $M$ is defined to be the graded module $M(n)$ with $M(n)_j = M_{n+j}$. The degree of a homogeneous element $f\in M$ will be denoted $|f|$.

An \textit{$A$-linear morphism of degree $j$} is a homomorphism $\alpha:M \to N$ of $A$-modules such that $\alpha(M_n) \subseteq N_{n+j}$ for each $n$. We let
	\[\grHom_A{(M,N)}_j
	\]
denote the $\bk$-space of all $A$-linear morphisms of degree $j$, and we set
	\[\grHom_A{(M,N)} = \bigoplus_{j\in \bbz} \grHom_A{(M,N)}_j.
	\]

The \textit{Hilbert series} of a graded $A$-module $M$ is defined to be the formal series
	\[\upH_M(t) = \sum_{j\in \bbz} \left(\dim_\bk{M_j}\right) t^j \in \bbz((t)).
	\]

Note that the ground field $\bk$ is both a left and right graded $A$-module, via the canonical projection $A \to A/A_+ \cong \bk$, where $A_+$ is the two-sided ideal of all elements of positive degree.

\subsection{Complexes and graded Ext}

Let $A$ be a connected $\bk$-algebra, and let $(X_\bt,\bd^X_\bt)$ and $(Y_\bt,\bd^Y_\bt)$ be chain complexes of graded $A$-modules. For each $i\in \bbz$ we set
	\[\grHom_A{(X,Y)}_i = \prod_{n\in \bbz} \grHom_A{(X_n,Y_{n+i})}.
	\]
If we then define
	\[\bd_i: \grHom_A{(X,Y)}_i \to \grHom_A{(X,Y)}_{i-1}
	\]
by
	\[\bd_i \left( (\alpha_n)_{n\in \bbz} \right) = \left( \bd^Y_{n+i}\circ \alpha_n - (-1)^i\alpha_{n-1} \circ \bd^X_n \right)_{n\in \bbz}
	\]
we obtain a chain complex $\left(\grHom_A{(X,Y)}_\bt,\bd_\bt\right)$ of graded vector spaces. Furthermore, this chain complex has a natural bigrading, defined by setting
	\[\grHom_A{(X,Y)}_{ij} = \prod_{n\in \bbz} \grHom_A{(X_n, Y_{n+i})}_j \quad \tn{for each $i,j\in \bbz$}.
	\]
Here $i$ represents the homological degree and $j$ the internal one. The subspaces of boundaries and cycles are homogeneous with respect to the internal degree and the differential preserves this degree, hence the cohomology
	\[\upH^\ast\left( \grHom_A{(X,Y)}\right)
	\]
inherits an internal degree.

The category of graded $A$-modules has enough projective and injectives. If $P_\bt$ is a projective resolution of a graded module $M$ and $I_\bt$ is an injective resolution of a graded module $N$, the standard argument (see, for example, \cite[Theorem 2.7.6]{Weibel1994}) can be adapted to show that there are isomorphisms
	\[\upH^\ast( \grHom_A{(M,I)}) \cong \upH^\ast ( \grHom_A{(P,I)}) \cong \upH^\ast( \grHom_A{(P,N)})
	\]
of bigraded $\bk$-spaces. Here $M$ and $N$ are viewed as complexes of graded $A$-modules concentrated in homological degree zero. The essential point is that $\grHom_A{(-,I)}$ and $\grHom_A{(P,-)}$ preserve quasi-isomorphisms, because $I$ is a bounded-above complex of injectives in the former case and $P$ is a bounded-below complex of projectives in the latter. The cohomology of $\grHom_A{(P,I)}$ is denoted $\grExt_A{(M,N)}$, and our isomorphisms above show that it can be computed by resolving either module. Note the bigrading:
	\[\grExt^i_A{(M,N)}_j = \upH^i\left( \grHom_A{(P,I)} \right)_j = \upH_{-i}\left( \grHom_A{(P,I)} \right)_j.
	\]

If $(X_\bt,\bd_\bt)$ is a chain complex of graded $A$-modules, the \textit{$m$-th homological shift} of $X$ is defined to be the chain complex $\Sigma^m X$ where $(\Sigma^m X)_i = X_{i-m}$ for all $i$. Its $i$-th differential
	\[(\Sigma^m X)_i \to (\Sigma^m X)_{i-1}
	\]
is defined
	\[x \mapsto (-1)^m \bd_{i-m}(x).
	\]
The \textit{$n$-th internal shift} of $X$ is defined to be the chain complex $X(n)$ where $X(n)_i = X_i(n)$ for all $i$.

\subsection{Dualization}

The linear dual of a vector space $X$ over a field $\bk$ will be denoted $X^\vee = \Hom_\bk{(X,\bk)}$. There is a $\bk$-linear map
	\[\Theta: X^\vee \otimes_\bk X^\vee \to (X \otimes_\bk X)^\vee, \quad \Theta( \varphi \otimes \psi) ( x \otimes y) = \varphi(x) \psi(y),
	\]
with $\varphi,\psi\in X^\vee$ and $x,y\in X$. If $X$ is finite-dimensional, then $\Theta$ is an isomorphism and we use it throughout the rest of this paper to identify the tensor product of duals with the dual of the tensor product.

\subsection{Minimal resolutions and presentations}

Every graded module $M$ over a connected $\bk$-algebra $A$ has a \textit{minimal} graded free resolution (see \cite[Chapter 1, Section 4]{PolishchukPositselski2005}). Furthermore, the bounded-below and locally finite-dimensional properties of $M$ pass to the free modules appearing in the resolution, and thus the dimensions
	\[\beta^A_{ij}(M) = \dim_\bk{\grExt_A^i{(M,k)}_{-j}},
	\]
called the \textit{Betti numbers} of $M$, are all finite. We define the \textit{bigraded Poincar\'e series} of $M$ to be the formal series
	\[\upP^A_M(s,t) = \sum_{i,j} \beta^A_{ij}(M) s^i t^j \in  \bbz[[s]]((t)).
	\]

A \textit{presentation} of $A$ is defined to be a choice of generating space $X$ of $A$ (as an algebra) and a choice of generating space $F$ of the kernel of the canonical surjection
	\[\varphi_X: \upT(X) \surj A
	\]
where $\upT(X)$ is the tensor algebra generated by $X$. The presentation is called \textit{minimal} if no proper subspace of $X$ generates the algebra, and no proper subspace of $F$ generates $\ker{\varphi_X}$. If $X$ has a basis $\{x_i\}_{i=1}^m$ and $F$ has a basis $\{f_j\}_{j=1}^n$, and if we want to display these bases explicitly, then we write
	\[A = \ang{x_1,\ldots,x_m \mid f_1,\ldots,f_n}.
	\]

Let $X$ be a finite-dimensional vector space over a field $\bk$ and suppose $f\in X \otimes X$. Let $\{x_1,\ldots,x_m\}$ be a basis of $X$ and expand $f$ as
	\[f = \sum_{1\leq i,j \leq m} \ell_{ij} (x_i \otimes x_j), \quad \ell_{ij} \in \bk.
	\]
Then the \textit{rank} of $f$ is defined to be the rank of the matrix $[\ell_{ij}]$. One can check easily that this definition does not depend on the chosen basis of $X$. The significance of this notion of rank is that it is linked to the structure of minimal free resolutions, as explained in the next proposition. But before stating this proposition, we lay out our convention which governs the matrix representations of boundary maps in free resolutions.

\begin{convention}\label{conv:matrix}
With respect to the standard bases, elements of the free right $A$-module $A_A^n$ will be represented as column vectors, and an $A$-linear map $A_A^n \to A_A^m$ will be represented as an $m\times n$ matrix which acts on the left.

On the other side, elements of the free left $A$-module $\ps{A}{A}^n$ will be represented as row vectors, and an $A$-linear map $\ps{A}{A}^n \to \ps{A}{A}^m$ will be represented as an $n\times m$ matrix which acts on the right. 

Free resolutions of right $A$-modules will be written with arrows pointing right-to-left, while free resolutions of left $A$-modules will be written with arrows pointing left-to-right. One reason we have adopted this convention is that it keeps the matrices in the order in which they would need to be multiplied if, for example, one was checking that the composition of two boundary maps in a resolution is zero.

According to this convention, dualizing a free resolution reverses the direction of the arrows, but does \textit{not} transpose the matrices.

This convention is essentially the one in \cite[Chapter I]{BellamyEtAl2016}.
\end{convention}

\begin{prop}\label{lem:hypersurfaceRes}
Suppose $A$ is a connected $\bk$-algebra with minimal presentation $A=\ang{x_1,\ldots,x_m \mid f}$ where $m\geq 2$, each $x_j$ has degree $1$, and $f$ has degree $2$. Write
	\[f = \sum_{j=1}^m x_j g_j \tn{ where $g_j = \sum_{i=1}^m \ell_{ij} x_i$ and $\ell_{ij} \in \bk$.}
	\]
\begin{enumerate}
\item If $f$ has rank $\geq 2$, then $A$ has global dimension $2$ and the minimal free resolution of $\bk_A$ is of the form
	\[0 \leftarrow A_A \xleftarrow{\begin{bmatrix}x_1 & \cdots & x_m \end{bmatrix}} A_A(-1)^m \xleftarrow{\begin{bmatrix}g_1 \\ \vdots \\ g_m\end{bmatrix}} A_A(-2) \leftarrow 0,
	\]
where we have observed Convention \ref{conv:matrix} in writing down the resolution.
\item The relation $f$ has maximal rank if and only if the minimal free resolution of $\ps{A}\bk$ is of the form
	\[0 \to \ps{A}A(-2) \xrightarrow{\begin{bmatrix}x_1 & \cdots & x_m \end{bmatrix}} \ps{A}A(-1)^m \xrightarrow{\begin{bmatrix}g_1 \\ \vdots \\ g_m\end{bmatrix}} \ps{A}A \to 0,
	\]
where we have observed Convention \ref{conv:matrix} in writing down the resolution.
\end{enumerate}

\begin{proof}
(1): We shall apply \cite[Theorem 5.3]{Dicks1985}, which states that to prove $A$ has global dimension $2$, it will suffice to show $f$ is not a square of a degree $1$ element. For this, suppose to the contrary that it is; say $f=g^2$ where $g = \sum_{i=1}^m \ell_i x_i$ for $\ell_i \in \bk$. Then
	\[f = \sum_{1\leq i,j\leq m} \ell_i \ell_j( x_i \otimes x_j).
	\]
But the matrix $[\ell_i \ell_j]$ has rank $1$, contradicting our assumption on the rank of $f$. Hence $f$ is not a square, and by the result quoted, this implies $A$ has global dimension $2$. The form of the free resolution then follows from, for example, \cite[Chapter I, Lemma 2.1.3]{BellamyEtAl2016}.

(2): If $f$ has maximal rank, then the set $\{g_1,\ldots,g_m\}$ is a basis of $A_1$ and the form the free resolution again follows from \cite[Chapter I, Lemma 2.1.3]{BellamyEtAl2016}. Conversely, if the displayed sequence is a resolution of $\ps{A}\bk$, then the set $\{g_1,\ldots,g_m\}$ is a basis for $A_1$, i.e., $f$ has maximal rank.
\end{proof}
\end{prop}

\subsection{Quadratic and Koszul algebras}\label{subsec:Koszul}

A connected $\bk$-algebra $A$ is called \textit{quadratic} if the canonical map
	\[\varphi_{A_1}: \upT(A_1) \to A
	\]
is surjective with kernel generated by the subspace $\ker{\varphi_{A_1}} \cap A_1^{\otimes 2}$. The \textit{quadratic part} of $A$ is defined to be the connected $\bk$-algebra with presentation
	\[A^\square = \ang{A_1 \mid \ker{\varphi} \cap A_1^{\otimes 2}};
	\]
thus the natural map $A^\square \to A$ is an isomorphism if and only if $A$ is quadratic.

If $A$ is quadratic and $\mu: A_1 \otimes_\bk A_1 \to A_2$ is the multiplication map, then the \textit{quadratic dual} of $A$ is defined to be the connected $\bk$-algebra with presentation
	\[A^! = \ang{A_1^\vee \mid \im{ \mu^\vee}}.
	\]
If $A$ is connected but not necessarily quadratic, then $A^!$ is defined to be the quadratic dual of $A^\square$, i.e., $(A^\square)^!$.

We say the algebra $A$ is \textit{Koszul} if it is quadratic and $\beta^A_{ij}(\bk)=0$ when $i\neq j$. The proofs of the following results can be found in \cite[Chapter 2]{PolishchukPositselski2005}.

\begin{prop}\label{prop:hilbertPoincare}
Let $A$ be a Koszul $\bk$-algebra.
\begin{enumerate}
\item There is an equality $\upH_A(t)  \upH_{A^!}(-t)=1$ of formal power series.
\item There is an equality $\upH_A(-st)  \upP^A_\bk(s,t)=1$ of formal power series.
\item The algebra $A$ has global dimension $n$ if and only if $A^!_n\neq 0$ and $A^!_{>n}=0$.
\end{enumerate}
\end{prop}

\subsection{Frobenius, Gorenstein, and AS-regular algebras}

Let $A$ be a connected algebra over a field $\bk$ and $M$ a graded left (right) $A$-module. We define
	\[D(M) = \grHom_k{(M,k)},
	\]
viewed as a graded right (left) $A$-module. Since our modules are locally finite-dimensional, the assignment $M \mapsto D(M)$ is a duality, i.e., there is a natural isomorphism $M \cong D^2(M)$ (given by evaluation).

The algebra $A$ is said to be \textit{Frobenius} if it is finite-dimensional and if there is an isomorphism
	\[A_A \cong D(\ps{A}A)(-n) \quad (n=\sup{A})
	\]
of graded right $A$-modules. Using the isomorphism $\ps{A}A \cong D^2(\ps{A}A)$, one can prove that this definition is left/right symmetric, i.e., we need not distinguish between \textit{left} and \textit{right} Frobenius algebras. This definition is a graded analog of a classical one; see, for example, \cite{MoorePeterson1973} and \cite{Smith1996} in the graded case, and \cite{CurtisReiner2006} for the classical theory.

The algebra $A$ is called \textit{Gorenstein} if the bigraded $\bk$-space $\grExt_A{(\bk_A,A_A)}$ is one-dimensional; in this case there are integers $d$ and $\ell$ such that
	\[\grExt^d_A{\left(\bk_A,\ps{A}A_A\right)} \cong \ps{A}\bk(\ell)
	\]
as graded left $A$-modules. Gorenstein algebras form an important class of algebras with an extensive literature, but we warn the reader that the definitions tend to vary --- our definition is from \cite[p. 24]{PolishchukPositselski2005}, which differs from the common one found in \cite{ArtinSchelter1987} (see also \cite[Chapter I, Definition 2.1]{BellamyEtAl2016}) by not requiring the algebra to have finite global dimension.

The following useful characterization of Frobenius and finite-dimensional Gorenstein algebras seems to be well-known, though a reference is difficult to locate in the form we require, so we supply a proof.

\begin{prop}\label{prop:socleFrob}
Let $A$ be a finite-dimensional connected $\bk$-algebra. The following are equivalent:
\begin{enumerate}
\item The algebra $A$ is Frobenius.
\item The algebra $A$ is Gorenstein.
\item The graded vector space
	\[\grSoc{A_A} = \{ a\in A : \tn{$a$ homogeneous and $a\cdot A_+=0$}\}
	\]
is one-dimensional.
\end{enumerate}
\begin{proof}
We first observe that in general there are isomorphisms
	\begin{equation}\label{eqn:socExt}
	\grSoc{A_A} \cong \grHom_A{(\bk_A,\ps{A}A_A)} \cong \grExt_A^0{(\bk_A,\ps{A}A_A)}
	\end{equation}
of graded left $A$-modules. Set $n=\sup{A}$.

(1) $\Rightarrow$ (2): By the Frobenius property, $(\otimes,\grHom)$-adjunction, and \eqref{eqn:socExt}, we have
\begin{align*}
\grSoc{A_A} &\cong \grHom_A{(\bk_A,A_A)} \\
&\cong \grHom_A{(\bk_A, D(\ps{A}A))}(-n) \\
&\cong \grHom_\bk{(\bk \otimes_A A,\bk)}(-n) \\
&\cong \bk(-n)
\end{align*}
as graded $\bk$-spaces. Thus $\grSoc{A_A}$ is a left ideal of $A$ concentrated in degree $n$, and therefore isomorphic to $\ps{A}\bk(-n)$ as a graded left $A$-module. To establish (2), we then appeal to \eqref{eqn:socExt} and the fact Frobenius algebras are self-injective.

(2) $\Rightarrow$ (3): The vector space $\grSoc A_A$ is nonzero since $A_n \subseteq \grSoc A_A$ for degree reasons. But then from \eqref{eqn:socExt} we conclude that this vector space is one-dimensional.

(3) $\Rightarrow$ (1): Again we must have $A_n \subseteq \grSoc A_A$, and since the latter space is one-dimensional, the containment is necessarily an equality. Then, since
	\[D(\grSoc A_A) \cong \frac{D(A_A)}{A_+D(A_A)},
	\]
Nakayama's lemma implies $D(A_A)$ is cyclic and generated in degree $-n$. Hence there is a surjection $\ps{A}A(n) \surj D(A_A)$, and since these two graded $\bk$-spaces have the same dimension, the surjection is necessarily an isomorphism.
\end{proof}
\end{prop}

Gorenstein algebras of finite global dimension can be characterized in terms of a symmetry condition on the minimal free resolution of $\bk$, as we now state. For a proof, see \cite[Chapter I, Section 2]{BellamyEtAl2016}.

\begin{prop}\label{prop:symmetricRes}
Let $A$ be a connected $\bk$-algebra of finite global dimension $d$ and let $F_\bt$ be the minimal graded free resolution of $\bk_A$. Then $A$ is a Gorenstein algebra if and only if there is an integer $\ell$ such that
	\[\Sigma^d \grHom_A{(F,\ps{A}A_A)}(-\ell)
	\]
is the minimal free resolution of $\ps{A}\bk$. In particular, if $A$ is a Gorenstein algebra, then the Betti numbers $\beta^A_{ij}(\bk)$ are symmetric in the sense that
	\[\beta^A_{i,j}(\bk) = \beta^A_{d-i,\ell-j}(\bk)
	\]
for all $i,j\in \bbz$.
\end{prop}

Combining the Gorenstein property with a finiteness and growth condition produces a class of algebras which has been intensely studied over the past three decades; the origins of this class are contained in \cite{ArtinSchelter1987}, \cite{ArtinTateVandenBergh1990}, and \cite{ArtinTateVandenBergh1991}, while a textbook account is given in \cite{BellamyEtAl2016}. We say a $\bk$-algebra is \textit{Artin-Schelter regular} or (\textit{AS-regular}) if it is Gorenstein, has finite global dimension, and has polynomial growth.

\subsection{Hilbert-Samuel multiplicity}\label{subsec:mult}

Let $(R,\mfm,k)$ be a commutative noetherian local ring. We write $R^\g$ for the associated graded ring $\gr_\mfm{R}$, and we write $R^\square$ in place of the cumbersome $(R^\g)^\square$.

If $M$ is a finitely-generated $R$-module of Krull dimension $d$, then we can write
	\[\upH_M(t) = \frac{h(t)}{(1-t)^d}
	\]
for some polynomial $h(t)\in \bbz[t]$ with $h(1)\neq 0$. We define the \textit{(Hilbert-Samuel) multiplicity of $M$} to be the integer
	\[e(M) = \begin{cases}
	h(1) & : \dim{M} = \dim{R}, \\
	0 & : \dim{M} < \dim{R}.
	\end{cases}
	\]
For further details and proofs see \cite[Chapter 4]{BrunsHerzog1998}.

The main classes of rings that this paper studies can now be defined. The \textit{codimension} of $R$, denoted $\codim{R}$, is defined to be the difference $\edim{R} - \dim{R}$.

\begin{itemize}
\item[(CM)] Let $R$ be a Cohen-Macaulay ring. Then Abhyankar shows in \cite{Abhyankar1967} that $e(R) \geq \codim{R} +1$. If equality holds, then $R$ is called a \textit{Cohen-Macaulay ring of minimal multiplicity}.
\item[(G)] Let $R$ be a Gorenstein ring with $e(R) \geq 3$. Then Sally shows in \cite[Corollary 3.2]{Sally1980-1} that $e(R) \geq \codim{R}+2$. If equality holds, then $R$ is called a \textit{Gorenstein ring of minimal multiplicity}.
\item[(CI)] Let $R$ be a complete intersection. Then $e(R) \geq 2^{\codim{R}}$ (see, for example, \cite[Theorem 14.10]{Matsumura1986}). If equality holds, then $R$ is called a \textit{complete intersection of minimal multiplicity}.
\end{itemize}

\section{Minimal multiplicity and $\Ext$-algebras}

Throughout this section $(R,\mfm,k)$ is a commutative noetherian local ring. We start by defining the $\Ext$-algebra $\sce_R$ of $R$, and we then investigate when the two algebras $\sce_R$ and $R^\square$ are quadratic duals of each other. This will lead to minimal multiplicity, and to our characterizations of this property via the structure of $\sce_R$.

\subsection{Ext-algebras of local rings}\label{subsec:extAlg}

If $F_\bt$ is a free resolution of $k$ with augmentation $\dev: F_\bt \to k$, then $\Hom_R{(F_\bt,\dev)}$ is a quasi-isomorphism inducing an isomorphism
	\[\upH^\ast(\Hom_R{(F,F)}) \xrightarrow{\cong} \upH^\ast( \Hom_R{(F,k)}) = \Ext_R{(k,k)}.
	\]
The composition product on $\bigoplus_{i\in \bbz}\Hom_R{(F,F)}_i$ descends to a product on the direct sum of its cohomology, and this structure is transported (via the isomorphism above) to $\bigoplus_{i\in \bbz} \Ext^i_R{(k,k)}$ making it a connected $k$-algebra. It is called the \textit{Ext-algebra} of $R$ and is denoted $\sce_R$. It is the universal enveloping algebra of a graded Lie algebra which will appear later in this paper.

\subsection{Quadratic duality and $\Ext$-algebras}

In \cite{Sjodin1976}, Sj\"odin explains how to compute the degree-$2$ relations of $\sce_R$. His computations begin by noting that $\sce_R$ is unchanged when we pass from $R$ to its $\mfm$-adic completion, and thus we may assume that we have an isomorphism $R \cong Q/\mfb$ where $(Q,\mfn,k)$ is a regular local ring and $\mfb$ is minimally generated by elements $f_1,\ldots,f_c$ contained in $\mfn^2$. If we suppose that $\mfn$ is minimally generated by $x_1,\ldots,x_n$, then for each $h=1,\ldots,c$ we can write
	\begin{equation}\label{eqn:relationsSjodin}
	f_h = \sum_{ 1\leq i \leq j \leq n} a_{h,ij} x_i x_j \quad \tn{with} \quad a_{h,ij}\in Q.
	\end{equation}
The cosets $x_i + \mfn^2$ form a basis of $\mfn/\mfn^2\cong \mfm/\mfm^2$, and if we set $\xi_i = (x_i+\mfn^2)^\vee$, then $\{\xi_1,\ldots,\xi_n\}$ is a basis of $\sce_R^1$. We then have

\begin{thm}[Sj\"odin {\cite[\S3]{Sjodin1976}}]\label{thm:sjodinPres}
Let the notation be as above, and let $\varphi: \upT(\sce_R^1) \to \sce_R$ be the natural map of connected $k$-algebras. The degree-$2$ homogeneous component of $\ker{\varphi}$ is generated by those elements of the form
	\[\sum_{1 \leq i < j \leq n} b_{ij}( \xi_i \otimes \xi_j + \xi_j \otimes \xi_i) + \sum_{i=1}^n b_{ii} \xi_i \otimes \xi_i \quad (b_{ij}\in k)
	\]
where the $\binom{n+1}{2} \times 1$ column vector $[b_{ij}]$ is contained in the kernel of the $c \times \binom{n+1}{2}$ matrix $[\cl{a}_{h,ij}]$ (here an overbar denotes an image in $k$).
\end{thm}

It follows that if $\sce_R$ is quadratic, then the relations Sj\"odin's theorem generate \textit{all} of the relations, and thus we have a complete presentation of $\sce_R$.

By pairing Sj\"odin's theorem along with an additional result in his paper \cite{Sjodin1976}, we obtain a proof of the next theorem, which is \textit{the} foundational theorem on which most of the oncoming results sit.

\begin{thm}\label{prop:dualOfExt}
Let $(R,\mfm,k)$ be a noetherian local ring. If $\sce_R$ is quadratic, then there is an isomorphism $\sce_R \cong (R^\square)^!$ of connected $k$-algebras.
\begin{proof}
We continue with the notation introduced above in the discussion leading to \eqref{eqn:relationsSjodin}. We write $\mfb^\ast$ for the ideal of $Q^\g$ generated by the initial forms of elements in $\mfb$, and we write $f^\ast$ for the initial form in $Q^\g$ of an element $f\in Q$.

Observe that since $\sce_R$ is quadratic, by \cite[Theorem 4]{Sjodin1976} the set $\{f_1^\ast,\ldots,f_c^\ast\}$ of initial forms in $Q^\g$ is linearly independent and contained in $\mfn^2/\mfn^3$; from this we deduce $\{f_1^\ast,\ldots,f_c^\ast\}$ is a basis of the vector space $\mfb^\ast \cap \mfn^2/\mfn^3$, which in turn means that we have an isomorphism
	\[R^\square \cong \frac{ Q^\g}{(f_1^\ast,\ldots,f_c^\ast)}
	\]
of graded $k$-algebras. If $X$ denotes the subspace of $Q^\g$ spanned by the initial forms $x_i^\ast$, then since $Q^\g$ is a polynomial algebra generated by $X$, we have $R^\square = \ang{ X \mid Y}$ where $Y$ is the subspace of $X\otimes X$ spanned by the elements $\sum_{ 1\leq i \leq j \leq n} \cl{a}_{h,ij} x_i^\ast \otimes x_j^\ast$ ($h=1,\ldots,c$) and $x_i^\ast \otimes x_j^\ast - x_j^\ast \otimes x_i^\ast$ ($1 \leq i < j \leq n$). We then have an exact sequence
	\[0 \to Y \xrightarrow{ \iota} X \otimes X \xrightarrow{ \mu} \frac{X \otimes X}{Y} \to 0
	\]
where $\mu$ denotes the multiplication map $R^\square_1 \otimes R^\square_1 \to R^\square_2$. Applying the $k$-linear dual yields the exact sequence
	\[0 \to \left(\frac{X \otimes X}{Y}\right)^\vee \xrightarrow{\mu^\vee} X^\vee \otimes X^\vee \xrightarrow{\iota^\vee} Y^\vee \to 0,
	\]
and since by definition $(R^\square)^! = \ang{X^\vee \mid \im{ \mu^\vee}}$, it follows that the degree-$2$ relations of $(R^\square)^!$ must be those elements of $X^\vee \otimes X^\vee$ which are in $\ker{\iota^\vee}$. But the elements of this kernel are exactly the elements of the form
	\[\sum_{1 \leq i < j \leq n} b_{ij}( (x_i^\ast)^\vee \otimes (x_j^\ast)^\vee + (x_j^\ast)^\vee \otimes (x_i^\ast)^\vee) + \sum_{i=1}^n b_{ii} (x_i^\ast)^\vee \otimes (x_i^\ast)^\vee \quad (b_{ij}\in k)
	\]
where the $\binom{n+1}{2} \times 1$ column vector $[b_{ij}]$ is in the kernel of the $c \times \binom{n+1}{2}$ matrix $[\cl{a}_{h,ij}]$. Thus, by Theorem \ref{thm:sjodinPres}, the algebras $\sce_R$ and $(R^\square)^!$ have the same presentation.
\end{proof}
\end{thm}

Thus from the quadratic condition on $\sce_R$ we deduce that this algebra is quadratic dual to $R^\square$; but what condition on $R$ will guarantee that $\sce_R$ is quadratic? An answer, at least in the artinian Gorenstein case, is minimal multiplicity. Indeed, we have the following stronger result:

\begin{thm}[Levin, Avramov {\cite[Theorem 3]{LevinAvramov1978}}]\label{thm:LevinAvramov}
Let $(R,\mfm,k)$ be an artinian Gorenstein local ring of embedding dimension $\geq 2$. If $R$ has minimal multiplicity (i.e., $\mfm^3=0$), then there is a presentation $\sce_R = \ang{\sce_R^1 \mid f}$ where $f$ is a degree-$2$ element of maximal rank.
\end{thm}

\begin{remark}
The maximal rank condition on the relation $f$ is not stated explicitly in Levin and Avramov's paper, so we outline how one deduces that $f$ has this property. Since $R$ is Gorenstein of minimal multiplicity, we have $\mfm^3=0$ and $\dim_k{\mfm^2}=1$; hence the multiplication map $\mu: \mfm/\mfm^2 \otimes \mfm/\mfm^2 \to \mfm^2$ can be viewed as a symmetric bilinear form, and the Gorenstein property is equivalent to nondegeneracy of this form. Hence $\im{\mu^\vee}$ is generated by an element $f \in (\mfm/\mfm^2)^\vee \otimes (\mfm/\mfm^2)^\vee$ of maximal rank. But since $\sce_R$ is quadratic, by Proposition \ref{prop:dualOfExt} there is an isomorphism $\sce_R \cong (R^\square)^!$ of connected $k$-algebras, and since
	\[(R^\square)^! = \ang{ (\mfm/\mfm^2)^\vee \mid \im{ \mu^\vee}},
	\]
the desired result follows.
\end{remark}

\subsection{Gorenstein rings of minimal multiplicity and $\sce_R$}
By coupling Theorem \ref{prop:dualOfExt} with Theorem \ref{thm:LevinAvramov}, we conclude that if $R$ is artinian Gorenstein of minimal multiplicity, then there is an isomorphism $\sce_R \cong (R^\square)^!$. Our Theorem 1 from the introduction shows that (with an added hypothesis on the global dimension of $\sce_R$) the existence of this isomorphism is also sufficient for the ring to have minimal multiplicity. We begin working toward a proof of this theorem by presenting some preliminary results on quadratic algebras with a single relation.

\begin{lem}[{\cite[Lemma 3.2]{Smith1996}}]\label{lem:poincareFrobenius}
Let $B$ be a connected algebra over a field $\bk$ with $\sup{B}=2$. Then $B$ is Frobenius if and only if $\dim_\bk{B_2}=1$ and the bilinear form
	\[\nu: B_1 \otimes B_1 \to B_2
	\]
given by multiplication is nondegenerate.
\end{lem}

\begin{lem}\label{lem:maximalNondegenerate}
Let $A$ be a quadratic algebra over a field $\bk$ with a single relation, say $A = \ang{A_1 \mid f}$. Then $f$ has maximal rank if and only if the bilinear form
	\[\nu: A_1^! \otimes A_1^! \to A_2^!
	\]
given by multiplication is nondegenerate.
\begin{proof}
Let $\{x_1,\ldots,x_m\}$ be a basis of $A_1$ and write
	\[f = \sum_{1 \leq i ,j \leq m} \ell_{ij}(x_i \otimes x_j), \ \ell_{ij}\in K.
	\]
The multiplication $\mu: A_1 \otimes A_1 \to A_2$ fits into an exact sequence
	\[0 \to \bk \cdot f \xrightarrow{\iota} A_1 \otimes A_1 \xrightarrow{\mu} \frac{A_1 \otimes A_1}{\bk \cdot f} \to 0
	\]
where $\iota$ is the inclusion. With respect to the bases $\{x_i^\vee \otimes x_j^\vee\}_{1\leq i,j \leq m}$ of $A_1^\vee \otimes A_1^\vee$ and $\{f^\vee\}$ of $\bk \cdot f^\vee$, the bilinear form $\iota^\vee:A_1^\vee \otimes A_1^\vee \to \bk \cdot f^\vee$ has matrix $[\ell_{ij}]$.

By definition $A^! = \ang{A_1^\vee\mid\im{\mu^\vee}}$, so the multiplication $\nu:A_1^! \otimes A_1^! \to A_2^!$ is equal to the canonical map
	\[A_1^\vee \otimes A_1^\vee \to \frac{ A_1^\vee \otimes A_1^\vee}{\im{\mu^\vee}}.
	\]
It fits into a commutative diagram
	\[\xymatrix{ A_1^\vee \otimes A_1^\vee \ar[dr]_{\iota^\vee} \ar[r]^\nu & \dfrac{A_1^\vee \otimes A_1^\vee}{\im{\mu^\vee}} \ar[d]^\cong \\
	&\bk \cdot f^\vee}
	\]
which shows $f$ has maximal rank if and only if $\iota^\vee$ is nondegenerate, if and only if $\nu$ is nondegenerate.
\end{proof}
\end{lem}

It is known that the Gorenstein and Frobenius properties are related through quadratic duality; see \cite[Section 9.3]{LuPalmieriWuZhang2008} and \cite[Proposition 5.10]{Smith1996}. We state here a version of this relation suited for our purposes.

\begin{prop}\label{prop:gradedGor}
Let $A$ be a connected algebra over a field $\bk$. The following are equivalent:
\begin{enumerate}
\item The algebra $A$ is quadratic Gorenstein of global dimension $2$.
\item The algebra $A$ is quadratic, and $A^!$ is Frobenius with $\sup{A^!}=2$.
\item The algebra $A$ has presentation $A = \ang{A_1 \mid f}$ where $f$ is a degree-$2$ element of maximal rank.
\end{enumerate}
\begin{proof}
Set $m = \dim_\bk{A_1}$.

(1) $\Leftrightarrow$ (3): In both cases (1) and (3), the minimal free resolution of $\bk_A$ (following Convention \ref{conv:matrix}) looks like
	\[0 \leftarrow A_A \xleftarrow{\bd_1} A_A(-1)^m \xleftarrow{\bd_2} A_A(-2) \leftarrow 0.
	\]
Indeed, in case (1) we use Proposition \ref{prop:symmetricRes} to compute the shifts, and in case (3) we use Proposition \ref{lem:hypersurfaceRes} to compute the global dimension. Hence in both cases we have $A = \ang{A_1 \mid f}$ where $f$ has degree $2$. If we write
	\[f = \sum_{j=1}^m x_j g_j \tn{ where $g_j = \sum_{i=1}^m \ell_{ij} x_i$ and $\ell_{ij} \in \bk$,}
	\]
then the differentials in the resolution are represented by the matrices
	\[\bd_1 = \begin{bmatrix} x_1 & \cdots & x_m \end{bmatrix}, \quad \bd_2 = \begin{bmatrix} g_1 \\ \vdots \\ g_m \end{bmatrix}.
	\]
Applying $\grHom_A{(-,A)}$ to the resolution produces a chain complex which is isomorphic to a homological and internal shift of
	\[0 \to \ps{A}A(-2) \xrightarrow{\begin{bmatrix}x_1 & \cdots & x_m \end{bmatrix}} \ps{A}A(-1)^m \xrightarrow{\begin{bmatrix}g_1 \\ \vdots \\ g_m\end{bmatrix}} \ps{A}A \to 0,
	\]
where we have again observed Convention \ref{conv:matrix} in writing the resolution left-to-right. Now, by Proposition \ref{lem:hypersurfaceRes}, the relation $f$ has maximal rank if and only if this complex resolves $\ps{A}\bk$, and by Proposition \ref{prop:symmetricRes} this occurs if and only if $A$ is Gorenstein.

(2) $\Leftrightarrow$ (3): If $A$ is \textit{any} quadratic $\bk$-algebra, then we have $\dim_\bk{A_2^!}=1$ if and only if $A$ has a single degree-$2$ relation. In view of this fact, the equivalence (2) $\Leftrightarrow$ (3) follows directly from Lemmas \ref{lem:poincareFrobenius} and \ref{lem:maximalNondegenerate}.
\end{proof}
\end{prop}

We now have enough to prove Theorem \ref{thm:1} from the introduction, restated here as

\begin{thm}\label{thm:localGor}
Let $(R,\mfm,k)$ be an artinian local ring of embedding dimension $\geq 2$. The following are equivalent:
\begin{enumerate}
\item The ring $R$ is Gorenstein of minimal multiplicity.
\item The ring $R$ is Gorenstein, the algebra $\sce_R$ has global dimension $2$, and there is an isomorphism $\sce_R \cong (R^\square)^!$ of connected algebras.
\item The algebra $\sce_R$ is quadratic Gorenstein of global dimension $2$.
\end{enumerate}

\begin{proof}

(1) $\Rightarrow$ (3): By Theorem \ref{thm:LevinAvramov}, the algebra $\sce_R$ has minimal presentation $\sce=\ang{\sce_R^1 \mid f}$ where $f$ is a degree-$2$ element of maximal rank. Thus, by Proposition \ref{prop:gradedGor}, the algebra $\sce_R$ is quadratic Gorenstein of global dimension $2$.

(3) $\Rightarrow$ (2): Since $\sce_R$ is 	quadratic, from Theorem \ref{prop:dualOfExt} we obtain the isomorphism $\sce_R \cong (R^\square)^!$, so all we need to prove is that $R$ is Gorenstein. For this, we note that Proposition \ref{prop:gradedGor} implies $R^\square$ is Frobenius with $\sup{R^\square}=2$. But then $R^\square$ and $R^\g$ coincide in degrees $\leq 2$, and so we must have $R^\square = R^\g$. Hence $R^\g$ is Frobenius, which means it is Gorenstein, and therefore $R$ is also Gorenstein.

(2) $\Rightarrow$ (1): Since $R$ is assumed Gorenstein, we need only prove $\mfm^3=0$. To do this, note that from the isomorphism $\sce_R \cong (R^\square)^!$ we may conclude the algebra $\sce_R$ is quadratic. But since it also has global dimension $2$, it is Koszul. Then applying quadratic duality to the isomorphism $\sce_R \cong (R^\square)^!$ yields $\sce_R^! \cong R^\square$, and hence $\sup{R^\square}=2$ by Proposition \ref{prop:hilbertPoincare}(3). Finally, since $\upH_{R^\square}(t) \geq \upH_R(t)$, we have that $\mfm^3=0$.
\end{proof}
\end{thm}

\begin{remarks}

Suppose $R$ is artinian with embedding dimension $\geq 2$. By the theorem and Proposition \ref{prop:gradedGor}, a quadratic presentation of the form $\sce_R = \ang{\sce_R^1 \mid f}$ with $f$ of maximal rank implies $R$ is Gorenstein of minimal multiplicity; we point this out because it is the converse of Levin and Avramov's Theorem \ref{thm:LevinAvramov}. However, the maximal rank condition on $f$ is essential for this implication to hold. Indeed, consider the ring
	\[R = \frac{k[x,y,z]}{(x^2,y^2,z^2,xy,xz)},
	\]
which decomposes as a fiber product $R \cong S\times_k T$ where $S=k[x]/(x^2)$ and $T = k[y,z]/(y^2,z^2)$. Since both $S$ and $T$ are complete intersections,  Sj\"odin \cite[Theorem 5]{Sjodin1976} shows that
	\[\sce_S \cong k\ang{\xi}  \quad \tn{and} \quad \sce_T \cong \frac{ k\ang{\upsilon,\zeta}}{(\upsilon\zeta +\zeta\upsilon)},
	\]
where $k\ang{\cdots}$ stands for the polynomial ring in \textit{non}commuting variables and each variable $\xi,\upsilon,\zeta$ has degree $1$. By \cite[Theorem 3.4]{Moore2009}, the $\Ext$-algebra of $R$ is then the free product (over $k$) of $\sce_S$ and $\sce_T$, namely
	\[\sce_R \cong \frac{k\ang{\xi,\upsilon,\zeta}}{(\upsilon\zeta +\zeta\upsilon)}.
	\]
Hence $\sce_R$ has a minimal presentation of the form $\ang{\sce_R^1 \mid f}$, but $R$ is not even Gorenstein.
\end{remarks}

\subsection{Complete intersections of minimal multiplicity and $\sce_R$}

With the Gorenstein case addressed, we now turn our attention to artinian complete intersections of minimal multiplicity and their $\Ext$-algebras.

Before stating the results, we elaborate briefly on an item mentioned in passing in Section \ref{subsec:extAlg}. The algebra $\sce_R$ (where $R$ is \textit{any} commutative noetherian local ring) is the graded universal enveloping algebra of a graded Lie algebra denoted $\pi^\ast(R)$ and called the \textit{homotopy Lie algebra} of $R$. We require several facts regarding the structure of $\pi^\ast(R)$ when $R$ is an artinian complete intersection of embedding dimension $m$. For these rings, $\pi^\ast(R)$ is concentrated in degrees $1$ and $2$ and $\pi^2(R)$ has dimension $m$ (see \cite[Corollary 7.1.5, Theorems 7.3.3, 10.2.1]{Avramov2010}). In particular, $\pi^2(R)$ is an abelian Lie subalgebra of $\pi^\ast(R)$ for degree reasons, and the universal enveloping algebra $\Gamma = \upU(\pi^2(R))$ is then a polynomial $k$-algebra generated by $m$ elements of degree $2$. The inclusion of Lie algebras $\pi^2(R) \inj \pi^\ast(R)$ induces an inclusion of associative algebras $\Gamma \inj \upU(\pi^\ast(R)) = \sce_R$, and the Poincar\'e-Birkhoff-Witt theorem (see \cite{Sjodin1980-2} for the version we require) then shows that $\sce_R$ is free as both a left and right graded $\Gamma$-module. Furthermore, for degree reasons the algebra $\Gamma$ is central in $\sce_R$, and by \cite[Theorem 5]{Sjodin1976} the quotient algebra $\Lambda = \sce_R/ (\Gamma_+ \cdot \sce_R)$ is an exterior algebra generated by $m$ elements of degree $1$.

We also require the following pair of results; statement (1) is due to Tate \cite[Theorem 6]{Tate1957}, while (2) is a result of Gulliksen \cite[Theorem 2.3]{Gulliksen1980}.

\begin{thm}[Tate, Gulliksen]\label{thm:CIFacts}
Let $(R,\mfm,k)$ be a noetherian local ring.
\begin{enumerate}
\item If $R$ is an artinian complete intersection, then $\upH_{\sce_R}(t) = 1/(1-t)^{\edim{R}}$.
\item The algebra $\sce_R$ has polynomial growth if and only if $R$ is a complete intersection.
\end{enumerate}
\end{thm}

We are now ready to prove Theorem \ref{thm:2} from the introduction, restated as

\begin{thm}\label{thm:CI}
Let $(R,\mfm,k)$ be an artinian local ring and set $m=\edim{R}$. The following are equivalent:
\begin{enumerate}
\item The ring $R$ is a complete intersection of minimal multiplicity.
\item The ring $R$ is a complete intersection, the algebra $\sce_R$ has global dimension $m$, and there is an isomorphism $\sce_R \cong (R^\square)^!$ of connected algebras.
\item[$(2')$] The ring $R$ is a complete intersection and the algebra $\sce_R$ has global dimension $m$.
\item The algebra $\sce_R$ is Koszul AS-regular of global dimension $m$.
\item[$(3')$] The algebra $\sce_R$ is Koszul and has polynomial growth.
\end{enumerate}
\begin{proof}
Set $\sce = \sce_R$.

(1) $\Rightarrow$ (2): As shown in the proof of \cite[Theorem 2.3]{Avramov1994}, the algebra $\sce$ is Koszul; the isomorphism $\sce \cong (R^\square)^!$ then follows from Theorem \ref{prop:dualOfExt}, and so all we have left is to compute the global dimension of $\sce$. For this, we note
	\[\upH_{R^\square}(t) = \frac{1}{\upH_{\sce}(-t)} = (1+t)^m
	\]
where the first equality is Proposition \ref{prop:hilbertPoincare}(1) and the second is from Theorem \ref{thm:CIFacts}. Hence $\sup{R^\square}=m$, which, by Proposition \ref{prop:hilbertPoincare}(3), means $\sce$ has global dimension $m$.

(2) $\Rightarrow$ (2$'$): Clear.

(2$'$) $\Rightarrow$ (3): According to Theorem \ref{thm:CIFacts}, the algebra $\sce$ has polynomial growth. Thus we need to show that $\sce$ is Gorenstein and Koszul.

To compute $\grExt_\sce{(k_\sce,\sce_\sce)}$ and show $\sce$ is Gorenstein, we note that we are in the situation described at the beginning of this section regarding the three algebras $\sce$, $\Gamma$, and $\Lambda$: The subalgebra $\Gamma$ of $\sce$ is a central polynomial subalgebra generated by $m$ elements of degree $2$, call them $x_1,\ldots,x_m$, and the quotient $\Lambda = \sce/(\Gamma_+ \cdot \sce)$ is an exterior algebra generated by $m$ elements of degree $1$. Furthermore, $\sce$ is free as both a right and left graded $\Gamma$-module.

Consider the first-quadrant cohomological change-of-rings spectral sequence associated to canonical morphism $\sce \surj \Lambda$:
	\[E_2^{pq} = \grExt_\Lambda^p{\left(k_\Lambda, \grExt_\sce^q{\left(\Lambda_\sce,\sce_\sce\right)}\right)} \Rightarrow \grExt^{p+q}_\sce{\left(k_\sce,\sce_\sce\right)}.
	\]
We've noted that $\Gamma$ is a polynomial algebra generated by $x_1,\ldots,x_m$, and hence the Koszul complex $K_\bt$ on the elements $x_1,\ldots,x_m$ resolves $k_\Gamma$. Since $\sce$ is flat as a left $\Gamma$-module, the complex $K \otimes_\Gamma\sce_\sce$ resolves $k\otimes_\Gamma \sce_\sce \cong \Lambda_\sce$, and we can therefore use $K \otimes_\Gamma\sce_\sce$ to compute $\grExt_\sce^q{\left(\Lambda_\sce,\sce_\sce\right)}$. We obtain
	\[\grExt_\sce^q{\left(\Lambda_\sce,\sce_\sce\right)} \cong \begin{cases}
	\Lambda_\Lambda(2m) & : q = m, \\
	0 & : q\neq m,
	\end{cases}
	\]
using the self-duality of $K$. Thus the spectral sequence collapses to the single row $q=m$, where we have
	\[E_2^{p,m} \cong \grExt_\Lambda^p{\left(k,\Lambda\right)}(2m).
	\]
But $\Lambda$ is finite-dimensional and Frobenius, so the only nonzero term on the second page of the spectral sequence is $E^{0,m}_2$, which is isomorphic to
	\[\grHom_\Lambda{(k,\Lambda)}(2m) \cong k(m)
	\]
since $\grHom_\Lambda{(k,\Lambda)} \cong k(-m)$. We conclude therefore that
	\[\grExt_\sce^i{(k_\sce,\sce_\sce)} \cong \begin{cases}
	k(m) & : i = m, \\
	0 & : i\neq m,
	\end{cases}
	\]
which means $\sce$ is Gorenstein.

To prove $\sce$ is Koszul, we apply \cite[Chapter 2, Theorem 2.5]{PolishchukPositselski2005} which states that to establish Koszulness of $\sce$ it is enough to show the Hilbert series $\upH_\sce(t)$ has degree $-m$. And this we already know: Since $R$ is an artinian complete intersection, Theorem \ref{thm:CIFacts} shows $\upH_\sce(t) = (1-t)^{-m}$.

(3) $\Rightarrow$ (3$'$): Clear.

(3$'$) $\Rightarrow$ (1): Polynomial growth implies $R$ is a complete intersection by Theorem \ref{thm:CIFacts}; thus we need only prove $R$ has multiplicity $2^m$. For this, we first note that from the isomorphism $\sce_R^! \cong R^\square$ in Theorem \ref{prop:dualOfExt} follows the equality $\upH_\sce(-t) \upH_{R^\square}(t)=1$ of Proposition \ref{prop:hilbertPoincare}(1). Then, since $R$ is a complete intersection, we have $\upH_\sce(t) = (1-t)^{-m}$ from Theorem \ref{thm:CIFacts}. The coefficient-wise inequality $\upH_{R^\square}(t) \geq \upH_R(t)$ thus turns into $(1+t)^m \geq \upH_R(t)$, from which it follows $e(R) \leq 2^m$. But we always have $2^m \leq e(R)$, and therefore $R$ has minimal multiplicity.
\end{proof}
\end{thm}

\section{Minimal multiplicity and linearity defect}

We now give a pair of characterizations of rings of minimal multiplicity in terms of a numerical invariant called linearity defect, and thereby answer in two special cases a question posed by Herzog and Iyengar in \cite{HerzogIyengar2005}.

Throughout this section $(R,\mfm,k)$ is a commutative noetherian local ring and $M$ is a finitely-generated $R$-module.

\subsection{Rational Poincar\'e series}

Some of the results in this section concern complete intersections and Golod rings, and it will be important to have at hand formulas for the Poincar\'e series of modules over them. We note that the Golod condition is defined in terms of a rational Poincar\'e series; indeed the ring $R$ is \textit{Golod} if there is an equation
	\[\upP^R_k(t) = \frac{(1+t)^{\edim{R}}}{1-\sum_{j=1}^{\codepth{R}}a_jt^{j+1}}
	\]
where the $a_j$'s are nonnegative integers and $\codepth{R} = \edim{R} - \depth{R}$. For more details about Golod rings, see \cite[Section 5]{Avramov2010}

In the next theorem, statement (1) is due to Tate \cite[Theorem 6]{Tate1957} and Gulliksen \cite[Corollary 4.2]{Gulliksen1974}, while the second statement is due to Ghione and Gulliksen \cite[Theorems 1 and 4]{GhioneGulliksen1975}.

\begin{thm}[Tate, Gulliksen, Ghione]\label{thm:ratlPoincare}
Let $(R,\mfm,k)$ be a noetherian local ring and $M$ a finitely-generated $R$-module.
\begin{enumerate}
\item If $R$ is a complete intersection, then
	\[\upP^R_k(t) = \frac{(1+t)^{\edim{R}}}{(1-t^2)^{\codim{R}}} \quad \tn{and} \quad \upP^R_M(t) = \frac{p_M(t)}{(1-t^2)^{\codim{R}}}
	\]
for some polynomial $p_M(t)\in \bbz[t]$.
\item If $R$ is Golod, then 
	\[\upP^R_{M}(t) = \frac{p_M(t)}{1-\sum_{j=1}^{\codepth{R}} a_j t^{j+1}}
	\]
for some polynomial $p_M(t)\in \bbz[t]$, and where the displayed denominator coincides with the one from $\upP^R_k(t)$.
\end{enumerate}
\end{thm}

\subsection{Linearity defect and a lemma}

Let $(F_\bt,\bd_\bt)$ be the minimal free resolution of $M$. Since $\bd(F_i) \subseteq \mfm F_{i-1}$ for each $i\geq 1$, the differentials fit into a subcomplex
	\[\scg^pF: \quad \cdots \to F_{p+1} \xrightarrow{\bd_{p+1}} F_p \xrightarrow{\bd_p} \mfm F_{p-1} \to \cdots \to \mfm^{p-1}F_1 \xrightarrow{\bd_1} \mfm^p F_0 \to 0
	\]
of $F$, one subcomplex for each $p\geq 0$. Then $\{ \scg^p F\}_{p\geq 0}$ is a descending filtration of $F$, and following Herzog and Iyengar \cite{HerzogIyengar2005} we call the associated graded complex the \textit{linear part} of $F$ and denoted it by $\lin^R{F}$. The number
	\[\ld_R{M} = \sup\{ i : \upH_i(\lin^R{F}) \neq 0\}
	\]
is called the \textit{linearity defect} of $M$. If $\ld_R{M} = 0$, then $M$ is called a \textit{Koszul module}. If $k$ is a Koszul module, then $R$ is called a \textit{Koszul ring}.

As shown in \cite[Proposition 1.5]{HerzogIyengar2005}, Koszulness of $R$ is equivalent to $R^\g$ being a Koszul algebra as defined in Section \ref{subsec:Koszul}, and in this case $\lin^R{F}$ provides a minimal free resolution of $k$ over $R^\g$. But then $\upP^{R^\g}_k(s,t) = \upP^R_k(st)$, and so by Proposition \ref{prop:hilbertPoincare}(2) we have
	\begin{equation}\label{eqn:frobergRelation}
	\upH_R(-t) \upP^R_k(t) = 1
	\end{equation}
when $R$ is a Koszul ring. Following Fitzgerald \cite{Fitzgerald1996}, we call this equality the \textit{Fr\"oberg relation}, and any ring that satisfies it is called a \textit{Fr\"oberg ring}. Hence the Koszul condition on $R$ implies the Fr\"oberg one, and it is an open question whether the converse holds. An affirmative answer is known in the case of complete intersections; see \cite[\S6.1, Proposition 6]{Conca2013}.

In \cite{Levin1985} (see also \cite[Theorem 6.3.6]{Avramov2010}), Levin discovered an analog of the Fr\"oberg relation which holds for \textit{any} $R$: He showed that for all $m\gg 0$ we have
	\begin{equation}\label{eqn:levinRelation}
	\upH^R_{\mfm^m}(-t) \upP^R_k(t) = \upP^R_{\mfm^m}(t).
	\end{equation}
Notice that this reduces to the Fr\"oberg relation when $m=0$. In \cite[Corollary 3.4]{Sega2001}, \c Sega gave a lower bound on the integer $m$ in terms of an invariant of $R$ called \textit{polynomial regularity}. To define this invariant, we present $R^\g$ as a quotient of a polynomial ring $k[x_1,\ldots,x_n]/\mfa$ where each $x_i$ has degree $1$ and $n = \edim{R}$. Then the \textit{polynomial regularity} of $R$ is defined to be the integer
	\[\polreg{R} = \reg_{k[x_1,\ldots,x_n]}{R^\g}
	\]
where the right-hand size is the Castelnuovo-Mumford regularity of $R^\g$ as a module over the polynomial ring $k[x_1,\ldots,x_n]$. \c Sega's result shows that if $m> \polreg{R}$, then Levin's relation \eqref{eqn:levinRelation} holds.

Before presenting the main technical lemma which will be used in the next section, we present a pair of results of \c Sega; they are proved in \cite[Proposition 3.2, Lemma 6.2]{Sega2013}.

\begin{prop}[\c Sega]\label{prop:sega}
Let $(R,\mfm,k)$ be a commutative noetherian local ring.
\begin{enumerate}
\item There is an inequality $\ld_R{\mfm^m} \leq \ld_R{k}$ for all $m$.
\item Let $R$ have Krull dimension $d$ and suppose $M$ is a finitely-generated $R$-module. Suppose that $\upP^R_M(t) = u(t)/g(t)$ for relatively prime polynomials $u(t),g(t)\in \bbz[t]$. If $M$ has finite linearity defect, then $g(-1)\neq 0$ and $e(M)/e(R) = u(-1)/g(-1)$.
\end{enumerate}
\end{prop}

We now present the main lemma. It is essentially a combination of the previous proposition and Levin's relation \eqref{eqn:levinRelation}. Its value is that it allows us to compute the multiplicities of certain rings from Poincar\'e series.

\begin{lem}\label{lem:finiteLinDefMult}
Suppose $(R,\mfm,k)$ is a commutative Cohen-Macaulay ring of dimension $d$ and that there exists integers $m,n$ and a polynomial $D(t)$ in $\bbz[t]$ such that $m>\polreg{R}$, $n\geq \dim{R}$,
	\[\upP^R_k(t)D(t) = (1+t)^n, \quad \tn{and} \quad \upP^R_{\mfm^m}(t)D(t) \in \bbz[t].
	\]
If $k$ has finite linearity defect, then
	\[D(t) = (1+t)^{n-d} g(t)
	\]
for some $g(t)$ in $\bbz[t]$ with $g(-1)=e(R)$.

\begin{proof}
Assume first that $d=0$. If we write $D(t) = (1+t)^q g(t)$ where $q$ is a nonnegative integer and $g(t)$ is a polynomial in $\bbz[t]$ with $g(-1) \neq 0$, then
	\[\upP^R_k(t) = \frac{ (1+t)^n}{g(t)(1+t)^q}.
	\]
An application of Proposition \ref{prop:sega}(2) with $M=k$ shows $n=q$ and $g(-1) =e(R)$, as desired.

Assume now $d>0$ and set $f(t) = \sum_{i=0}^{m-1} \dim_k(\mfm^i/\mfm^{i+1})t^i$ and $p(t) = \upP^R_{\mfm^m}(t) D(t)$. Then
	\begin{equation}\label{eq:first}
	\upH_R(t) = f(t) + t^m \upH^R_{\mfm^m}(t) = f(t) + t^m \frac{\upP^R_{\mfm^m}(-t)}{\upP^R_k(-t)} = \frac{f(t)(1-t)^n+p(-t)t^m}{(1-t)^n}
	\end{equation}
where the second equality is a consequence of \eqref{eqn:levinRelation}. However, we also have $\upH_R(t) = h(t)/(1-t)^d$ for some $h(t) \in \bbz[t]$ with $e(R) = h(1)$, and equating the right-hand side of \eqref{eq:first} with $h(t)/(1-t)^d$ results in
	\begin{equation}\label{eq:second}
	p(-t)t^m = (1-t)^{n-d} \left( h(t) - f(t) (1-t)^d\right).
	\end{equation}
This implies $(1-t)^{n-d}$ must divide $p(-t)$; say $p(-t) = (1-t)^{n-d}u(-t)$ for some $u(t) \in \bbz[t]$. Replacing $p(-t)$ in \eqref{eq:second} with this factorization and canceling gives
	\[u(-t) t^m  = h(t) -  f(t)(1-t)^d.
	\]
Since $d>0$ and $e(R) = h(1)$, evaluating both sides of the equality at $t=1$ shows $e(R) = u(-1)$.

As in the dimension $0$ case, write $D(t) = (1+t)^q g(t)$ where $q\geq 0$ and $g(t)$ is a polynomial in $\bbz[t]$ with $g(-1)\neq 0$. Then
	\[\upP^R_{\mfm^m}(t) = \frac{p(t)}{D(t)} = \frac{(1+t)^{n-d}u(t)}{(1+t)^q g(t)}.
	\]
Now, the ideal $\mfm^m$ has the same multiplicity as $\mfm$ and by Proposition \ref{prop:sega}(1) it has finite linearity defect; thus from Proposition \ref{prop:sega}(2) we conclude $q=n-d$ and that $g(-1) = u(-1)$. But $u(-1)=e(R)$ from the previous paragraph, so we are done.
\end{proof}
\end{lem}

\subsection{Characterizations of minimal multiplicity via linearity defect}

We begin with complete intersections. For these rings, \c Sega showed the Koszul property of $R$ and minimal multiplicity are equivalent, and that if $R^\g$ is Cohen-Macaulay, then these two properties are also equivalent to $k$ having finite linearity defect \cite[Theorem 6.5]{Sega2013}. Also, as we noted when we defined Fr\"oberg rings, it is known that the Koszul and Fr\"oberg conditions are equivalent for complete intersections. Thus our contribution in the next theorem is to prove \c Sega's Cohen-Macaulay hypothesis on $R^\g$ is not necessary, and  Lemma \ref{lem:finiteLinDefMult} will make short work of it.

\begin{thm}\label{thm:CILinDef}
Let $(R,\mfm,k)$ be a complete intersection. The following are equivalent:
\begin{enumerate}
\item The ring $R$ is Koszul.
\item The ring $R$ is Fr\"oberg.
\item The ring $R$ is a complete intersection of minimal multiplicity.
\item The $R$-module $k$ has finite linearity defect.
\end{enumerate}

\begin{proof}
In view of the discussion preceding the theorem, we need only prove (4) $\Rightarrow$ (3). From Theorem \ref{thm:ratlPoincare} we have
	\[\upP^R_k(t)(1-t^2)^{\codim{R}} = (1+t)^{\edim{R}} \quad \tn{and} \quad\upP^R_{\mfm^m}(t) (1-t^2)^{\codim{R}}  \in \bbz[t].
	\]
for all $m\geq 0$. Now, in the notation of Lemma \ref{lem:finiteLinDefMult} we have $n=\edim{R}$ and $D(t) = (1-t^2)^{\codim{R}}$; then the conclusion of the lemma is that there is $g(t)\in \bbz[t]$ with $g(t)(1+t)^{\codim{R}} = D(t)$ and $e(R) = g(-1)$. But our only choice is $g(t) =(1-t)^{\codim{R}}$, and so $e(R) = g(-1) = 2^{\codim{R}}$.
\end{proof}
\end{thm}

We now address rings with are Cohen-Macaulay and Golod. As with the complete intersection case, we again find ourselves building on work of \c Sega; in \cite[Theorem 7.2]{Sega2013} she proved the equivalences (1) $\Leftrightarrow$ (3) $\Leftrightarrow$ (4) in the next theorem, but only in the presence of a Cohen-Macaulay hypothesis on the associated graded ring. We use Lemma \ref{lem:finiteLinDefMult} to move this Cohen-Macaulay hypothesis from the associated graded ring to the ring itself. The equivalence (1) $\Leftrightarrow$ (2) is also new.

\begin{thm}\label{thm:CMGolodLinDef}
Let $(R,\mfm,k)$ be Cohen-Maculay and Golod. The following are equivalent:
\begin{enumerate}
\item The ring $R$ is Koszul.
\item The ring $R$ is Fr\"oberg.
\item The ring $R$ is Cohen-Macaulay of minimal multiplicity.
\item The $R$-module $k$ has finite linearity defect.
\end{enumerate}
\begin{proof}

(1) $\Rightarrow$ (2): We have already mentioned that the Koszul property implies the Fr\"oberg one.

(1) $\Rightarrow$ (4): Clear.

(3) $\Rightarrow$ (1): Minimal multiplicity guarantees that the associated graded ring $R^\g$ is Cohen-Macaulay; see \cite[Theorem 2]{Sally1977}. Now apply \c Sega's result mentioned before this theorem.

The implications that are left are (4) $\Rightarrow$ (3) and (2) $\Rightarrow$ (3). For them, set $d=\dim{R}$, $n=\edim{R}$, $c=\codim{R}$. From Theorem \ref{thm:ratlPoincare} we have
	\begin{equation}\label{eqn:k}
	\upP^R_k(t) \left(1- \sum_{j=1}^{c} a_j t^{j+1}\right) = (1+t)^n,
	\end{equation}
and
	\begin{equation}\label{eqn:m}
	\upP^R_{\mfm^m}(t) \left(1- \sum_{j=1}^{c} a_j t^{j+1}\right)\in \bbz[t]
	\end{equation}
for some nonnegative integers $a_j$ and for all $m$.

(4) $\Rightarrow$ (3): We apply Lemma \ref{lem:finiteLinDefMult}. Comparing its statement with \eqref{eqn:k} and \eqref{eqn:m} above shows $D(t) = 1- \sum_{j=1}^{c} a_j t^{j+1}$ and that $D(t) = g(t) (1+t)^c$ for some $g(t) \in \bbz[t]$ with $e(R) = g(-1)$. By equating coefficients we determine $g(t) = 1-ct$, and so $e(R) = c+1$.

(2) $\Rightarrow$ (3): Combining \eqref{eqn:k} with the Fr\"oberg relation \eqref{eqn:frobergRelation} yields
	\[\upH_R(t) = \frac{1+\sum_{j=1}^{c} (-1)^{j}a_jt^{j+1}}{ (1-t)^n}.
	\]
Now, if $\upH_R(t) = h(t)/(1-t)^{d}$ for $h(t) \in \bbz[t]$ with $e(R) = h(1)$, then
	\[h(t)(1-t)^{c} = 1+\sum_{j=1}^{c} (-1)^ja_jt^{j+1},
	\]
and then equating coefficients shows $h(t) = 1 +ct$. Hence $e(R) = c+1$.
\end{proof}
\end{thm}

\section*{Acknowledgements}

This paper is based off of the author's Ph.D. thesis, completed under the supervision of L. L. Avramov at the University of Nebraska-Lincoln. The author warmly thanks Professor Avramov for his guidance and encouragement during the completion of this work, and also thanks the anonymous referee who offered valuable suggestions to improve the paper.

This work was partly supported by NSF grant DMS-1103176.

\bibliographystyle{amsplain}
\bibliography{bib}

\end{document}